\documentclass{siamart220329}
\usepackage{amssymb}
\usepackage{array}
\makeatletter
\@ifundefined{definition}{\newsiamthm{definition}{Definition}}{}
\makeatother
\usepackage{booktabs}
\usepackage{graphicx}
\usepackage{tikz}
\usetikzlibrary{arrows.meta,positioning,fit,backgrounds,calc}
\usepackage{algorithm}
\usepackage{algpseudocode}
\usepackage{enumitem} 
\usepackage[]{mdframed}
\usepackage[skins,breakable]{tcolorbox}
\graphicspath{{figures/}}
\usepackage{multirow}
\usepackage{utfsym}
\usepackage{tablefootnote}
\usepackage{caption}

\title{DBDSVDMR3:\ An MRRR-Based Bidiagonal SVD Solver Built from DSTEGR\thanks{Submitted to the editors October 2026.
\funding{R.S.\ is supported by NSF grant MSPRF 2402027.}}}
\author{James Demmel\thanks{University of California Berkeley} \and Inderjit Dhillon\thanks{University of Texas at Austin and Google} \and Sai Surya Duvvuri\thanks{University of Texas at Austin and Core Automation} \and Julien Langou\thanks{University of Colorado Denver} \and Osni Marques\thanks{Lawrence Berkeley National Laboratory} \and Ryan Schneider\footnotemark[2]}

\headers{DBDSVDMR3: An MRRR Bidiagonal SVD Solver}{Demmel, Dhillon, Duvvuri, Langou, Marques, and Schneider}

\begin{document}
\maketitle

\begin{abstract}
The optimal algorithm for the symmetric tridiagonal eigenvalue problem is Multiple Relatively Robust Representations (MRRR or MR$^3$), which is capable of computing $k$ eigenvalue/eigenvector pairs in only $O(kn)$ floating-point operations. While MRRR is nearly 30 years old, it has never successfully been adapted to the bidiagonal SVD. This is not for lack of trying; two previous attempts produced implementations that were considered for LAPACK but eventually abandoned, a result of the fact that neither was reliable enough when run on test problems. The question therefore remained:\ is there a problem with the theory of MRRR for the bidiagonal SVD or were there simply bugs in the previous implementations that went undetected? This paper resolves this question by presenting \texttt{DBDSVDMR3} -- a new implementation of the MRRR-based algorithm for the bidiagonal SVD originally proposed by Willems and Lang [\textit{ETNA} 2012], which, for the first time, performs well on a variety of test problems while achieving significant speedups  over both LAPACK's bisection-based solver and, in certain cases, divide-and-conquer for the bidiagonal SVD 
(thereby indicating that the theory is sound). \texttt{DBDSVDMR3} will be available in the next release of LAPACK. 
In this paper, we discuss the design of our  implementation, which was obtained with the help of Claude Code, and exhibit its performance on a diverse set of test problems. 
\end{abstract}



\section{Introduction}\label{sec: intro}
Multiple Relatively Robust Representations (MRRR or MR$^3$) is the so-called ``holy-grail" algorithm for the symmetric tridiagonal eigenvalue problem. A special version of bisection + inverse iteration,\footnote{The technical details of the algorithm are deferred to Section 2.} MRRR can compute $k$ eigenpairs of an $n \times n$ symmetric tridiagonal $T$ in only $O(kn)$ floating-point operations (flops), which is optimal. While the algorithm is somewhat technical, its central idea is remarkably simple:\ if the eigenvalues of $T$ are computed accurately enough -- and if additionally \textit{relative} eigenvalue gaps are large -- then the eigenvectors obtained via inverse iteration will be mutually orthogonal by default, forgoing the need to run expensive, manual (re)orthogonalization. \\
\indent MRRR was first introduced in the 1997 PhD thesis of Dhillon \cite{Dhillon1997} and subsequently presented in a pair of papers by Dhillon and Parlett \cite{DhillonParlett2004,DhillonParlett2004b}. It was initially implemented in  LAPACK as \texttt{DSTEGR} \cite{DhillonParlettVoemel2006}, though this routine has since been superseded by \texttt{DSTEMR}. In the nearly three decades since its introduction, MRRR has cemented itself as part of the fastest algorithm for the symmetric eigenvalue problem, which first reduces a symmetric matrix $A$ to the tridiagonal $T$ (via Householder similarity transformations) and subsequently diagonalizes $T$ with MRRR. Accordingly, the algorithm has been implemented in a wide range of software libraries, including ScaLAPACK \cite{LAWN168} and MAGMA (see \texttt{sy/hegvr}). \\
\indent Following the work of Dhillon and Parlett, attention turned to using MRRR to compute the singular value decomposition (SVD) of a (real) bidiagonal matrix $B$. Ideally, an MRRR-based algorithm would, in $O(kn)$ flops for $B$ an  $n \times n$ bidiagonal, compute $k$ singular value/vector triplets $(\sigma_i,v_i,u_i)$ such that 
\vspace{1mm}
\begin{enumerate}[itemsep=.5em]
    \item $\|Bv_i - \sigma_i u_i\|_2,  \|u_i^TB-\sigma_iv_i^T \| =  O(\varepsilon n \|B\|)$
    \item $\left| v_i^Tv_j\right|, \left| u_i^Tu_j\right| = O(\varepsilon n)$ if $i \neq j$
    \item $\left| 1-\|v_i\|^2 \right|,  \left|1-\|u_i\|^2 \right| = O(\varepsilon n)$,
\end{enumerate}
\vspace{2mm}
where $\varepsilon$ is the machine precision used and $\|\cdot \|$ is an appropriate matrix norm. The first item here guarantees a small residual error while the second two imply (floating-point) orthogonality for the singular vectors. Note that we do not require high relative accuracy for the computed singular values $\sigma_i$, which is attainable for nonsingular bidiagonals via certain algorithms (e.g., the QR iteration  \cite{DemmelKahan1990}). \\
\indent We should expect to be able to use MRRR here  given the natural connection between the bidiagonal SVD and the symmetric tridiagonal eigenvalue problem. In particular, a bidiagonal $B$ is associated with three symmetric tridiagonal matrices:\ (1) $B^TB$, (2) $BB^T$, and $(3)$ the $2n \times 2n$ Golub-Kahan matrix
\begin{equation}\label{eqn: gk_definition}
    T_{\text{GK}} = P^T \begin{pmatrix} 0 & B \\ B^T & 0 \end{pmatrix} P,
\end{equation}
where $P$ is a \textit{perfect-shuffle} permutation (see \cite{GolubKahan1965}). Note that the diagonal of $T_{\text{GK}}$ is zero. If $B = U\Sigma V^T$ is an SVD, then 
\begin{equation}
    B^TB = V\Sigma^2V^T, \; \; BB^T = U\Sigma^2U^T, \; \; \text{and} \; \; \; T_{\text{GK}} = (PQ) \begin{pmatrix} -\Sigma & 0 \\ 0 & \Sigma \end{pmatrix} (PQ)^T, 
\end{equation} 
where we have
\begin{equation}\label{eqn: GK_eigenvectors}
    Q = \frac{1}{\sqrt{2}} \begin{pmatrix} U & U \\ -V &  V \end{pmatrix}.
\end{equation}
In other words, the SVD of $B$ can be obtained by diagonalizing either $B^TB$ and $BB^T$ or $T_{\text{GK}}$, possibly via MRRR. Similar to the tridiagonal case, this could be used in a fast algorithm for the general SVD, which first reduces an arbitrary matrix to bidiagonal form \cite{Barlow_reduction} and subsequently computes its SVD via MRRR. \\
\indent Unfortunately, such an algorithm has never successfully materialized. Indeed, the LAPACK routine for standard (bidiagonal) bisection + inverse iteration  -- i.e., \texttt{DBDSVDX} -- is still ``...to be replaced with a version of the Multiple Relative Robust Representation algorithm" according to its documentation.\footnote{See \href{https://www.netlib.org/lapack/explore-html/dc/d73/group__bdsvdx.html}{https://www.netlib.org/lapack/explore-html/dc/d73/group\_\_bdsvdx.html}}\\
\indent Historically, it has been considered ill-advised to compute the bidiagonal SVD by applying MRRR as a ``black box" to one of its associated tridiagonals. In the case of approach one, which computes $V$ and $U$ by independently diagonalizing $B^TB$ and $BB^T$, the resulting singular vectors may be insufficiently ``coupled," meaning the residual error $\|Bv_i-\sigma_i u_i\|_2$ can be large. On the flip side, extracting singular vectors from the eigenvector matrix of $T_{\text{GK}}$ risks orthogonality; that is, $U^TU$ and $V^TV$ may not be close to the identity, even if $Q$ is orthogonal to working precision. These pitfalls were highlighted in work of Großer and Lang \cite{GrosserLang2003,GrosserLang2005}, which extended MRRR to the bidiagonal SVD via a coupled approach, working with all three of $B^TB$, $BB^T$, and $T_{\text{GK}}$ simultaneously. An implementation of their algorithm was proposed for LAPACK in 2005 \cite{LAWN166} but was ultimately abandoned (without much explanation but likely due to failure cases). \\
\indent Willems and Lang again took up the question of MRRR for the bidiagonal SVD in a series of papers beginning in 2010 \cite{WillemsLang2012,WillemsLang2012ETNA,WillemsLang2014,WillemsLang2011Blocks}, including the Householder-prize-winning PhD thesis of Willems \cite{Willems2010}. This work not only reformulated the proof of MRRR for the symmetric tridiagonal eigenvalue problem but also introduced MR3-GK, an algorithm for the bidiagonal SVD that could (provably, under some mild assumptions) attain the guarantees listed above. For the latter, their main insight was the following:\ in contrast to traditional lore, MRRR \textit{can} be successfully applied as a black box\footnote{In fact, some care must be taken to ensure success. See Section 3.} to the Golub-Kahan matrix. This resulted in a new implementation that was once again submitted for addition to LAPACK. Unfortunately, it was also eventually abandoned; as noted in work of Marques, Demmel, and Vasconcelos \cite{MDV2020}, it can fail on even some simple examples. \\
\indent This left open the following question, which was included in a list of open problems in numerical linear algebra created as part of a semester program at the Simons Institute \cite[Section 3.4]{Simons2026}.
\vspace{2mm}
\begin{mdframed}
\noindent \textbf{Key Question:} Is there a fundamental flaw with using MRRR to compute the SVD of a bidiagonal matrix -- in particular via the MR3-GK approach -- or have the previous implementations failed because of bugs in their code?
\end{mdframed}
\vspace{2mm}
The difficulty in answering this question lies primarily in how technical MRRR is as an algorithm and hence how complicated corresponding implementations are. \\
\indent In this paper, we resolve this question by presenting a new implementation of the MRRR-based algorithm of Willems and Lang (i.e., MR3-GK). This is the first version of MRRR for the bidiagonal SVD that passes the LAPACK test code. Moreover, it performs well on a collection of ``torture tests" maintained by Marques \cite{STCollection},  attaining error ratios (with respect to backward error and loss of orthogonality, in units of $\varepsilon n$) below 100 on 110 of the 125 matrices tested (88\%), with median ratios about ten times those of the existing LAPACK routines for the bidiagonal SVD; the remaining problems feature  bidiagonals with adversarially clustered eigenvalues, for which MRRR loses accuracy --  e.g., it returns a decomposition, but with at least one error ratio above 100 (see Section~\ref{sec:stc}). 
On these same test problems, \texttt{DBDSDC} and \texttt{DBDSQR} are accurate throughout  while \texttt{DBDSVDX} fails twice. \\
\indent We obtain our implementation by allowing Claude Code  \cite{ClaudeCode} to build a version of the MR3-GK algorithm out of the original \texttt{DSTEGR} code for the symmetric tridiagonal eigenvalue problem. Specifically, our starting point is the \texttt{DSTEGR} of LAPACK 3.0 (1999) and its kernels \texttt{DLARRV}, \texttt{DLAR1V}, \texttt{DLARRB}, and \texttt{DLARRF}, not the current \texttt{DSTEGR}, which since LAPACK 3.1 is a wrapper around the substantially different \texttt{DSTEMR} \cite{DhillonParlettVoemel2006}. In this way, our work also serves as a one-algorithm case study into the use of automated code generation tools to produce new and/or modified software for numerical linear algebra. The resulting implementation has been merged into the reference LAPACK repository as \texttt{DBDSVDMR3} (pull request 1361, October 2026) and will be included in the next LAPACK release.  \\
\indent Table~\ref{tab: algorithm_comp} compares our version of MRRR with the other algorithms for the bidiagonal SVD present in LAPACK. Each of these satisfies the three guarantees mentioned above (and may also compute high relative accuracy singular values).  The algorithms also differ in what determines the accuracy of their computed singular \emph{vectors}:\ for \texttt{DBDSDC} and \texttt{DBDSVDX}, the error in $u_i$ and $v_i$ is bounded in terms of the \emph{absolute} gap $\min_{j\ne i}|\sigma_i-\sigma_j|$, whereas for zero-shift QR it is bounded in terms of the \emph{relative} gap $\min_{j\ne i}|\sigma_i-\sigma_j|/\sigma_i$, as conjectured by Demmel and Kahan \cite{DemmelKahan1990} and proved by Deift, Demmel, Li, and Tomei \cite{DDLT1991}. MRRR on $T_{\mathrm{GK}}$ belongs in principle to the second class, since its eigenvectors are determined by relative gaps of the representation they are computed from. In \texttt{DBDSVDMR3}, this holds within each block produced by the preprocessing discussed in Section~\ref{sec:prep}, but the absolute splitting criterion used to obtain these blocks reduces the overall guarantee to one dependent on the absolute gap. Historically, divide-and-conquer \cite{GuEisenstat1995} (i.e., \texttt{DBDSDC}) is considered the fastest of these algorithms if a full SVD is requested. 

\noindent \begin{minipage}{13cm}
\renewcommand{\arraystretch}{1.2}
    \centering
    \captionof{table}{Algorithms for the bidiagonal SVD. Complexity bounds are for the problem of computing $k$ singular value/vector triplets for $B \in {\mathbb C}^{n \times n}$.}
    \label{tab: algorithm_comp}
    \resizebox{\linewidth}{!}{
    \begin{tabular}{lcccc}   
            \toprule
             \multirow{2}{*}{Algorithm} & LAPACK & Relatively accurate  & 
             Singular value gap & Worst-case \\
             & routine & singular values? 
             & determining vector accuracy & arithmetic complexity \\
            \midrule 
            QR Iteration \cite{DemmelKahan1990} & \texttt{DBDSQR} & \checkmark 
            & relative  & $O(n^3)$ \\
            Divide-and-Conquer \cite{GuEisenstat1995} & \texttt{DBDSDC} & \scalebox{0.8}{\usym{2613}} 
            & absolute & $O(n^3)$ \\
            Bisection + Inverse Iteration & \texttt{DBDSVDX} & \checkmark 
            & absolute & $O(nk^2)$ \\
            \midrule
            \multirow{2}{*}{MRRR}  & \texttt{DBDSVDMR3} & \multirow{2}{*}{\hspace{.5mm} \scalebox{0.8}{\usym{2613}}\footnote{It is possible to formulate MRRR to get high accuracy singular values (with respect to relative error), but our current implementation does not do this. See Section 3 for the details.}} 
            & \multirow{2}{*}{ \hspace{.1mm} absolute 
            } & \multirow{2}{*}{$O(nk)$\footnote{ \texttt{DBDSVDMR3} does not yet attain this complexity, as it computes $n$ singular vectors even when a subset is requested and therefore always performs $O(n^2)$ operations. Future work will provide the option to return only $k$ singular value/vector triplets at the lower $O(kn)$ cost (see Section~\ref{sec:blas}).}}  \\
            & (2026) & 
            & & \\
            \bottomrule
        \end{tabular}
        }
\renewcommand{\arraystretch}{1}
\end{minipage}

\vspace{3mm}

\indent The contributions of the remaining sections can be summarized as follows. In Section~\ref{sec:bg}, we present the details of the original MRRR algorithm before, in Section~\ref{sec:adapt}, summarizing the modifications necessary to obtain a successful algorithm for the bidiagonal SVD from \texttt{DSTEGR}. In Section~\ref{sec:acc}, we present the results of our implementation on a diverse test suite of bidiagonals, demonstrating that it is both accurate and competitive in speed with existing options for the bidiagonal SVD (i.e., those listed in Table~\ref{tab: algorithm_comp}). We also examine how the speed comparison depends on the underlying BLAS used. Finally, in Section~\ref{sec:ai} we discuss our experience using AI tools to develop numerical linear algebra software and present recommendations for other researchers. 
\section{Background} \label{sec:bg}
We begin by reviewing two MRRR algorithms:\ (1) the original version for the symmetric tridiagonal eigenvalue problem introduced by Dhillon and Parlett \cite{Dhillon1997,DhillonParlett2004,DhillonParlett2004b} and (2) MR3-GK, the algorithm of Willems and Lang for the bidiagonal SVD \cite{WillemsLang2012ETNA,Willems2010}. Both are presented at a high-level; for the technical details, we refer readers to the listed references.

\subsection{MRRR for symmetric tridiagonals} The original MRRR algorithm, presented here as Algorithm~\ref{alg:mr3}, computes $k$ eigenpairs of a symmetric tridiagonal $T$. This is done by traversing a \textit{representation tree} and computing each eigenpair (via bisection and inverse iteration) from a  \textit{relatively robust representation} (RRR) of $T$, which determines the pair to high relative accuracy. In this context, a representation of $T$ can consist of any of the following:
\begin{enumerate}
    \item The entries of $T$ itself.
    \item The entries of an $LDL^T$ or $UDU^T$ factorization of $T$, where $L$ is lower bidiagonal, $U$ is upper bidiagonal, and $D$ is diagonal. 
    \item The entries of a \textit{twisted} factorization, obtained by fusing parts of $LDL^T$ and $UDU^T$ factorizations of $T$.
\end{enumerate}
A representation is \textit{relatively robust} with respect to the eigenpair $(\lambda_i,v_i)$ if applying small, relative perturbations to its entries results  in only small, relative changes to $\lambda_i$ and  (normwise) to $v_i$. Any of the aforementioned representations can constitute an RRR for a subset of the eigenpairs $(\lambda_i,v_i)$ (see \cite{DhillonParlett2004}), though we note that the tridiagonal $T$ itself is usually not sufficient.\footnote{See the appropriately titled chapter ``Tridiagonals are inadequate" from the PhD thesis of Dhillon \cite{Dhillon1997}. } 

\begin{algorithm}[t]
\caption{High-level MRRR for the tridiagonal eigenproblem
(after \cite[Alg.~2.1]{WillemsLang2012ETNA}) \\
\textbf{Input:} (1) $T \in {\mathbb C}^{n \times n}$ a symmetric tridiagonal. \\
\phantom{Input2 }  (2) $I_0 \subset [n]$ a set of indices corresponding to the desired eigenvalues. \\
\phantom{Input3 } (3) A pre-set gap tolerance $\textit{gaptol} > 0$. \\
\textbf{Output:} Eigenpairs $(\lambda_i,v_i)$ of $T$ with $i \in I_0$.
} \label{alg:mr3}
\begin{algorithmic}[1]
\State Compute a root representation $M_0$ for $T$ (after shifting if necessary)
\State S = $(M_0, I_0, 0)$ \Comment{Initialize representation tree}
\While{$S \neq \emptyset$}
  \State Remove $(M,I,\tau)$ from $S$
    \State Compute local eigenvalues $\lambda_i^{\text{loc}}$ of $M$ with $i \in I$ and classify $I = I_1 \cup \cdots \cup I_m$ \indent  into singletons and clusters according to $\textit{gaptol}$
    \For{$j = 1:m$}
    \If{$I_j = \left\{ i \right\}$} \Comment{Case 1 --  Singleton}
    \State Refine $\lambda^{\text{loc}}_i$ and compute a corresponding eigenvector $v_i$ via inverse 
    \indent \indent \; \; iteration (applied to $M$)
    \State $\lambda_i = \lambda^{\text{loc}}_i + \tau$
    \Else \Comment{Case 2 -- Cluster}
    \State Choose a shift $\tau'$ near the cluster 
    \State Compute a relatively robust representation $M'$ of $M-\tau'$
    \State Add $(M',I_j,\tau+\tau')$ to $S$
    \EndIf
  \EndFor
\EndWhile \\
\Return $(\lambda_i,v_i)$ for $i \in I_0$.
\end{algorithmic}
\end{algorithm}

\indent The key result underpinning MRRR is the following (for a proof of this proposition, see \cite{Parlett_book}).
\begin{proposition}\label{prop: gap_result}
    Let $T$ be a symmetric tridiagonal with eigenvalues $\lambda_1, \ldots, \lambda_n$ and corresponding eigenvectors $v_1, \ldots, v_n$. If $(\lambda,v)$ is an approximation to the eigenpair $(\lambda_i,v_i)$ of $T$, with $\|v\|_2=1$, and $\theta_i$ is the angle between $v$ and $v_i$, then
    \[ \sin(\theta_i) \leq \frac{\|Tv-\lambda v \|_2}{\min_{j \neq i} |\lambda-\lambda_j|}. \]
\end{proposition}

\indent Proposition~\ref{prop: gap_result} implies that, if we can compute an approximation $(\lambda,v)$ to $(\lambda_i,v_i)$ such that $\|Tv-\lambda v\|_2 \leq O(\varepsilon n|\lambda_i|)$, which is possible from a corresponding RRR by definition, and if additionally $\lambda_i$ is \textit{relatively} well separated from the other eigenvalues of $T$, then $v$ will be nearly aligned with the associated true eigenvector $v_i$. Hence, eigenvectors obtained from relatively robust representations of $T$ are orthogonal by default. \\
\indent The only catch here is that we need relative eigenvalue gaps to be sufficiently large, but this is easy to manage; if a cluster of eigenvalues is detected, we can simply shift near it to increase relative gaps. This describes the main loop of MRRR:\ given an RRR and a set of eigenvalues to compute, the algorithm determines which eigenvalues are ``singletons" and which belong to clusters -- according to a pre-set relative gap tolerance -- and computes eigenvectors for the singletons while shifting near each cluster and computing a new, corresponding RRR for it.  \\
\indent Naturally, this requires subroutines for (1) quickly computing the eigenvalues of an RRR and (2) computing a new RRR for $T - \tau I$ given an RRR for $T$. These can be handled with fast, differential qd algorithms \cite{dqd,WillemsLang2011Blocks}. We again omit the details of these subroutines for brevity (and since they will not be necessary to understand the extension to the bidiagonal setting). \\
\indent The notation used in Algorithm~\ref{alg:mr3} follows Willems and Lang \cite{WillemsLang2012ETNA}.  In particular, $(M,I,\tau)$ is a triple satisfying the following:\ $M$ is a representation of $T - \tau I$ that is an RRR for eigenvalues with indices in $I \subseteq [n]$, which we refer to as the \textit{local} eigenvalues corresponding to $M$. The \emph{root representation} $M_0$ corresponds to the first of these triples in the representation tree.  

\subsection{MR3-GK for bidiagonals} 
\begin{algorithm}[t]
\caption{MR3-GK for the bidiagonal SVD
\cite[Alg.~4.1]{WillemsLang2012ETNA} \\
\textbf{Input:} (1) $B \in {\mathbb R}^{n \times n}$ an upper bidiagonal matrix. \\
\phantom{Input2 }  (2) $I_0 \subset [n]$ a set of indices corresponding to the desired singular values. \\
\textbf{Output:} Singular value/vector triplets  $(\sigma_i,v_i,u_i)$ of $B$ with $i \in I_0$.
} \label{alg:mr3-gk}
\begin{algorithmic}[1]
\State Form the Golub-Kahan matrix $T_{\text{GK}}$ corresponding to $B$.
\State Run MRRR (i.e., Algorithm~\ref{alg:mr3}) on $T_{\text{GK}}$, using its entries as the root representation and requesting only eigenpairs with positive eigenvalues. 
\State For each eigenvector $q_i$ of $T_{\text{GK}}$, corresponding to eigenvalue $\sigma_i$ with $i \in I_0$, set
\[ \begin{pmatrix} u_i \\ v_i \end{pmatrix} = \sqrt{2} P^Tq_i\]
for $P$ a $2n \times 2n$ perfect-shuffle permutation matrix. \\
\Return $(\sigma_i,v_i,u_i)$ for $i \in I_0$.
\end{algorithmic}
\end{algorithm}

As mentioned in Section~\ref{sec: intro}, it was long believed that computing the SVD of a bidiagonal $B$ by applying Algorithm~\ref{alg:mr3} as a black box to either $BB^T$, $B^TB$, or the corresponding Golub-Kahan matrix $T_{\text{GK}}$ was fundamentally flawed. In \cite{WillemsLang2012ETNA}, Willems and Lang proved this to be unfounded, introducing MR3-GK (shown here as Algorithm~\ref{alg:mr3-gk}), an algorithm that computes singular value/vector triplets of $B$ by running MRRR on $T_{\text{GK}}$. \\
\indent Under some basic assumptions (see the next subsection), Willems and Lang demonstrate that success for this approach is guaranteed as long as the following hold:
\begin{enumerate}
    \item The entries of $T_{\text{GK}}$ are used as the root representation in MRRR.
    \item Any representation $M$ used in the algorithm is a small, relative entry-wise perturbation away from a representation that corresponds to a tridiagonal with constant diagonal. Willems and Lang refer to this as the ``nearly constant diagonal" condition. 
\end{enumerate}
The motivation for item two is the following:\ any symmetric tridiagonal with constant diagonal is simply a shifted Golub-Kahan matrix; accordingly,  its invariant subspaces have ``GK structure," meaning each eigenvector has block form matching \eqref{eqn: GK_eigenvectors} (see \cite[Lemma 4.4]{WillemsLang2012ETNA}). Hence, the invariant subspaces obtained from a representation $M$ satisfying the nearly-constant-diagonal condition have \textit{near}-GK structure, and it is safe to read off left/right singular vectors as done in line 3 of Algorithm~\ref{alg:mr3-gk} with minimal impact to their orthogonality. This resolves the main roadblock to running MRRR on $T_{\text{GK}}$. \\
\indent In practice, we will apply some preprocessing to $T_{\text{GK}}$ (e.g., zeroing out small entries and breaking the tridiagonal into independent blocks) in between lines 1 and 2 in Algorithm~\ref{alg:mr3-gk}. The details of this preprocessing -- and its implementation  -- are discussed in Section~\ref{sec:adapt}. For simplicity, we assume for the remainder of the paper that $B$ is an upper bidiagonal with real entries.\footnote{In the general SVD application, we can always choose to reduce arbitrary $A \in {\mathbb C}^{n \times n}$ to a real, upper bidiagonal matrix.}

\subsection{Unified theory of correctness} 
In addition to MR3-GK, the other major contribution of Willems and Lang was a reworking of the proof of correctness for the original MRRR algorithm. In \cite{WillemsLang2012}, they guarantee success for Algorithm~\ref{alg:mr3} under five basic assumptions (see Table~\ref{tab:conform}). Since Algorithm~\ref{alg:mr3-gk} calls Algorithm~\ref{alg:mr3} with minor modifications, this theory extends to it as well. Indeed, the same five assumptions -- plus the nearly constant diagonal condition mentioned above -- will guarantee correctness for MR3-GK. We state this as an informal theorem.

\begin{theorem}[informal version of {\cite[Thm.~4.5]{WillemsLang2012ETNA}}]\label{thm:wl45}
If the conditions listed in Table~\ref{tab:conform} are satisfied, then the singular vectors computed by Algorithm~\ref{alg:mr3-gk} are numerically
orthogonal and the coupled residuals
$\|Bv_i - u_i \sigma_i\|_2$,
$\|u_i^TB- \sigma_iv_i^T\|_2$ are small $($on the order of
$O(\varepsilon n \|B\|_2/\textit{gaptol})$ for $\varepsilon$ the machine precision$)$.
\end{theorem}

\begin{table}[!t]\caption{Willems and Lang's requirements for the correctness of Algorithms~\ref{alg:mr3} and Algorithm~\ref{alg:mr3-gk}.  See \cite{WillemsLang2012ETNA} for a precise statement of each requirement.}
\label{tab:conform}
\centering
\begin{tabular}{lp{9.5cm}}
\toprule
Requirement & Meaning \\
\midrule
\textsc{rrr} & Every representation used is an RRR. \\
\textsc{elg} & Element growth is (conditionally) bounded. \\
\textsc{relgaps} & Eigenvalue partitions are made, with respect to relative gaps, according to \textit{gaptol} (in line 5 of Algorithm~\ref{alg:mr3}). \\
\textsc{shiftrel} & Shifts are computed with mixed relative stability. \\
\textsc{getvec} & Eigenvectors are computed accurately from nodes of the representation tree. \\
\textsc{gk structure} & Invariant subspaces have near-GK structure
(bidiagonal only). \\
\bottomrule\end{tabular}\end{table}

\section{Adapting DSTEGR to the Bidiagonal SVD}
\label{sec:adapt} Our routine \texttt{DBDSVDMR3} implements Algorithm~\ref{alg:mr3-gk} on top of \texttt{DSTEGR}, the original LAPACK version of MRRR for symmetric tridiagonals (which, as mentioned in the introduction, is taken from LAPACK 3.0). 
This is not the first implementation of Algorithm~\ref{alg:mr3-gk}; as discussed above, Willems and Lang supplied an initial version, which was never adopted due to poor performance, specifically with respect to accuracy -- see \cite{MDV2020}. \\
\indent The purpose of this section is to discuss how to adapt \texttt{DSTEGR} for use on bidiagonals. The main modifications are described in detail below. Each of these changes was implemented by an AI coding assistant working under the direction of the authors. See Section~\ref{sec:ai} for a more detailed discussion of this process.  
\begin{enumerate}
    \item Preprocessing:\ Before calling MRRR, \texttt{DBDSVDMR3} breaks the Golub-Kahan matrix $T_{\text{GK}}$ into independent subproblems (where possible). We can increase the number of these subproblems by first zeroing out small off-diagonal entries of $T_{\text{GK}}$,  though this may come at the expense of only guaranteeing high absolute (as opposed to relative) accuracy in the resulting singular values.
    \item Using $T_{\text{GK}}$ as the root representation:\ The main insight of Willems and Lang is that MRRR can be run as a black box on the Golub-Kahan matrix $T_{\text{GK}}$ provided it is used as the root representation. Hence, we modify \texttt{DSTEGR} to do exactly this if it detects that the input matrix has GK-structure (e.g., if it has a constant diagonal). 
    \item Respecting the GK-structure of $T_{\text{GK}}$'s eigenvectors:\ The original \texttt{DSTEGR} code uses a simple trick to speed up the computation of eigenvectors -- if it detects that a vector's entries are decaying, it truncates and sets the remainder to zero. This was determined to be too aggressive for the Golub-Kahan matrix  (whose eigenvectors interleave corresponding singular vectors) and is therefore disabled.
\end{enumerate}
We show in the next section that \texttt{DBDSVDMR3} -- incorporating these changes -- performs well on a diverse set of test problems. These results indicate that the theory put forward by Willems and Lang is sound (thereby resolving our key question).\\
\indent Throughout this section, $B$ is an $n\times n$ (real) upper bidiagonal matrix with diagonal $d \in {\mathbb R}^n$, superdiagonal $e \in {\mathbb R}^{n-1}$, and singular values $\sigma_1\ge\dots\ge\sigma_n$. We assume a standard floating-point setting, where $\varepsilon$ is the machine precision used and $\mathrm{sfmin}$ is the smallest normalized floating-point number. If IEEE double precision is used, these take the values $\varepsilon=2^{-52}\approx2.2\times10^{-16}$ and $\mathrm{sfmin}\approx2.2\times10^{-308}$. Finally, $\| \cdot \|$ denotes the maximum entry norm (to be differentiated from the spectral norm $\|\cdot \|_2$, equal to the largest singular value). For reference, a flow of \texttt{DBDSVDMR3} is presented in Figure~\ref{fig:flow}. 


\subsection{Preprocessing}
\label{sec:prep}
On input, \texttt{DBDSVDMR3} forms $T_{\text{GK}}$ and processes it according to the following steps. Willems and Lang propose a similar approach to preprocessing, albeit in terms of the corresponding bidiagonal $B$, in \cite{WillemsLang2012ETNA}. 
\begin{enumerate}
    \item First, $B$ (and hence $T_{\text{GK}}$) is scaled so that $\max_{i}\{|d_i|,|e_i|\}$ falls within the interval $\left[ \sqrt{\text{sfmin}},1/\sqrt{\text{sfmin}} \right]$. This is the scaling inherited from \texttt{DSTEGR}, and its purpose is to keep the squares and products of entries formed by dqds and by the twisted factorizations from over/underflowing. It does not by itself zero out small entries, which are dealt with in the next step. 
    \item Next, entries of $d$ and $e$, which are interleaved on the off-diagonals of $T_{\text{GK}}$, are set to zero if they satisfy $|d_i|,|e_i| \leq  \varepsilon\|B\|$. Since this effectively applies an absolute perturbation to $B$ of size $O(\varepsilon \|B\|)$ it eliminates the possibility of getting high relative accuracy singular values from \texttt{DBDSVDMR3} (and can be skipped if desired).
    \item The upshot of zeroing out entries is that we can now easily split $T_{\text{GK}}$ into independent subproblems. An example -- where zeroing out the second diagonal entry of a $4 \times 4$ bidiagonal naturally decomposes the corresponding $8 \times 8$ Golub-Kahan matrix into a direct sum of smaller tridiagonals -- is shown below:
  \begin{equation*}
       \resizebox{.9\textwidth}{!}{$T_{\text{GK}} = \begin{pmatrix}
             0& d_1 & & & &  & & \\
            d_1& 0 & e_1 & & & & &  \\
            & e_1& 0& d_2 & &  & & \\
            & &  d_2& 0 & e_2 &  & & \\
            & & & e_2 & 0 & d_3 & & \\
            & & & & d_3 & 0  & e_3 &\\
            & & & & & e_3 &0 & d_4 \\
            & & & & & & d_4  & 0
        \end{pmatrix} \xrightarrow{d_2 = 0} \begin{pmatrix}
              {\color{blue} 0 }&  {\color{blue} d_1} & & & & & & \\
             {\color{blue} d_1} &  {\color{blue} 0} &  {\color{blue} e_1} & & &  & & \\
            & {\color{blue}  e_1} &  {\color{blue} 0} &  & &  & & \\
            & &  & {\color{red} 0} & {\color{red} e_2} & & & \\
            & & & {\color{red} e_2} & {\color{red} 0} & {\color{red} d_3} & & \\
            & & & & {\color{red} d_3} & {\color{red} 0} & {\color{red} e_3} & \\
            & & & & &  {\color{red} e_3} & {\color{red} 0} & {\color{red} d_4} \\
            & & & & & & {\color{red} d_4} & {\color{red} 0}  .        
        \end{pmatrix} = \begin{pmatrix}  {\color{blue} T_1} & \\ & {\color{red} T_2} \end{pmatrix}$.}
    \end{equation*}
    \item From this example, we see that it is possible to get a Golub-Kahan matrix of odd dimension after splitting (and in fact this always occurs when a diagonal entry of $B$ is set to zero). While we have assumed to this point that Golub-Kahan matrices have even dimension, this is not necessary for the theory. Indeed, if such a matrix has dimension $2k + 1$ then it corresponds to \eqref{eqn: gk_definition} for $B$ a $k \times (k+1)$ bidiagonal. Importantly, an odd-dimensional Golub-Kahan matrix is still an RRR for its spectrum and therefore can be used as a root representation in \texttt{DBDSVDMR3}.\footnote{In fact, Demmel and Kahan prove that bidiagonal matrices determine their singular values to high relative accuracy by showing the same of Golub-Kahan matrices, of any dimension, and their eigenvalues \cite[Theorem 2]{DemmelKahan1990}.}
    In this case, GK-structure looks like the following:\ if $B = U\Sigma V^T$ is a \textit{reduced} SVD with $U \in {\mathbb C}^{k \times k}$, $\Sigma \in {\mathbb R}^{k \times k}$, and $V \in {\mathbb C}^{(k+1)\times k}$, then
    \begin{equation}\label{eqn: odd_GK_structure}
        T_{\text{GK}} = (PQ) \begin{pmatrix} -\Sigma & & \\ & 0 & \\ & & \Sigma \end{pmatrix} (PQ)^T \quad \text{for} \quad \; Q = \frac{1}{\sqrt{2}}  \begin{pmatrix} U & 0 &U \\ -V& z& V\end{pmatrix}, 
    \end{equation}
    where $P$ is a $(2k+1) \times  (2k+1)$ perfect-shuffle permutation and $z \in {\mathbb C}^{k+1}$ belongs to the null space of $B$ (note that $T_{\text{GK}}$ will always have a zero eigenvalue if its dimension is odd). 
    \item Some care must be taken when handling these zero eigenvalues. In particular, recall that the singular values of $B$ repeat in $\pm$ pairs as eigenvalues of $T_{\text{GK}}$. Hence, while -- in the above example -- we obtain a zero eigenvalue from each of $T_1$ and $T_2$, the corresponding bidiagonal $B$ has only one zero singular value. Moreover, it is clear from \eqref{eqn: odd_GK_structure} that we cannot get both a corresponding left and right singular vector from an odd-sized Golub-Kahan matrix. Accordingly, we must compute a zero eigenvalue/eigenvector pair from both $T_1$ and $T_2$ to construct the corresponding left/right singular vectors of $B$. To be specific (again for our $4 \times 4$ example with $d_2 = 0$), $B$ is a direct sum of a $1 \times 2$ upper bidiagonal and a $3 \times 2$ lower bidiagonal -- i.e., 
    \begin{equation*}
        B = \begin{pmatrix}
             {\color{blue} d_1} &  {\color{blue} e_1} & & \\
            & & {\color{red} e_2} &  \\
            & & {\color{red} d_3} & {\color{red} e_3} \\
            & & & {\color{red} d_4} 
        \end{pmatrix} = \begin{pmatrix}  {\color{blue} B_1} & \\ & {\color{red} B_2}\end{pmatrix}.
    \end{equation*}
    $T_1$ and $T_2$ are then, respectively, the Golub-Kahan matrices corresponding to $B_1$ and $B_2^T$, and their null vectors contribute $z_1$ and $z_2$ satisfying $B_1z_1 = 0$ and $B_2^Tz_2 = 0$. Padding these vectors with zeros yields the right singular vector $(z_1 \; 0 )^T$ and the left singular vector $(0 \; z_2)^T$ of $B$ corresponding to its zero singular value.
\end{enumerate}

We emphasize that the loss of relative accuracy in \texttt{DBDSVDMR3} is a consequence of the splitting step above and not the singular value computation itself. Whenever no entry of $B$ falls below $\varepsilon\|B\|$, this step is skipped and dqds runs on the whole matrix. As a result, every singular value is delivered to full relative accuracy. A simple example makes the contrast with divide-and-conquer, the other algorithm for the bidiagonal SVD that does not guarantee high relative accuracy, clear. Suppose that $d_i=1$ and $e_i=100$ for all $i$. Since the product of the singular values is $\det(B)=1$, the smallest singular value must be approximately $10^{-2(n-1)}$. For $n=40$, \texttt{DBDSVDMR3} and \texttt{DBDSQR} both return $\sigma_{\min}=9.999\times10^{-79}$, and for $n=100$ both return $1.0\times10^{-199}$, agreeing to all printed digits, whereas \texttt{DBDSDC} (with vectors, so that divide-and-conquer is used) returns $2.6\times10^{-27}$ in both cases, correct only to the absolute level $O(\varepsilon\|B\|)$ of its backward error. With this in mind, we note that we plan to update \texttt{DBDSVDMR3} in the future to allow users to request higher accuracy in the singular values, if desired, in which case the splitting step discussed above is either handled more delicately or is skipped altogether.

\begin{figure}[!ht]\centering
\definecolor{flowblue}{RGB}{220,233,245}\definecolor{flowgreen}{RGB}{230,242,224}\definecolor{flowgreendk}{RGB}{50,120,60}
\definecolor{floworange}{RGB}{253,227,201}\definecolor{flowyellow}{RGB}{253,245,204}
\begin{tikzpicture}[font=\small, node distance=4mm,
  box/.style={draw=black!80, rounded corners=2pt, line width=0.6pt, align=center, inner sep=4pt, text=black},
  pre/.style={box, fill=flowblue, text width=0.84\linewidth},
  mr3/.style={box, fill=flowgreen, text width=0.80\linewidth},
  half/.style={box, fill=flowgreen, text width=0.375\linewidth},
  root/.style={box, fill=floworange, text width=0.80\linewidth},
  post/.style={box, fill=flowyellow, text width=0.84\linewidth},
  inbox/.style={box, fill=black!8, text width=0.70\linewidth, font=\small\bfseries},
  flow/.style={draw=black, -{Triangle[length=2mm,width=2.2mm]}, line width=0.8pt}]
\node[inbox, text width = 9cm] (in) {Input: Real, upper bidiagonal $B$ and an index set $I$.};
\node[pre, below=of in] (pre) {\textsc{preprocessing}\\ Form the Golub-Kahan matrix $T_{\text{GK}}$. \\ Scale if needed and set negligible entries to zero.};
\node[font=\small\bfseries, text=flowgreendk, below=7mm of pre] (mr3title) {Run MRRR on each independent block of $T_{\text{GK}}$ (via modified \texttt{DSTEGR}) };
\node[root, below=2mm of mr3title] (root) {Use the block  itself as the root representation. \\ Request only eigenpairs corresponding to nonnegative eigenvalues.};
\node[mr3, below=of root] (bis) {Compute (or refine) relevant eigenvalues. \\
Group into clusters and singletons according to \textit{gaptol}.};
\node[half, below left=4mm and -0.4\linewidth of bis] (iso) {\textsc{singleton}\\ Compute a corresponding eigenvector via inverse iteration (applied to a twisted factorization).};
\node[half, below right=4mm and -0.4\linewidth of bis] (clu) {\textsc{cluster}\\ Shift near the cluster, compute a new RRR, and return to the previous step.};
\node[post, below=12mm of iso.south -| bis] (extract) {Read off left/right singular vectors according to \eqref{eqn: GK_eigenvectors}. \\ Singular vectors corresponding to zero singular values can be obtained \\
from independent (odd-sized) blocks.};
\node[post, below=of extract] (out) {Sort and  undo scaling (if used). \\
If exceptions occur or an error is flagged, default to \texttt{DBDSDC}};
\node[inbox, below=of out, text width = 11cm] (final) {Output: Singular value/vector triplets $(\sigma_i, v_i,u_i)$ of $B$ with $i \in I$. };
\newlength{\mrboxextend}
\setlength{\mrboxextend}{.5cm}  

\coordinate (mr3boxbottom) at ($(clu.south)+(0,-\mrboxextend)$);

\begin{scope}[on background layer]
\node[draw=flowgreendk, dashed, line width=1pt, rounded corners=4pt, fill=flowgreen!40,
      inner sep=5pt, minimum width=12cm,
      fit=(mr3title)(root)(bis)(iso)(clu)(mr3boxbottom)] (mr3box) {};
\end{scope}

\draw[flow] (in) -- (pre); 
\draw[flow] (pre) -- (mr3box.north); 
\draw[flow] (root) -- (bis); 
\draw[flow] (bis) -- (iso); \draw[flow] (bis) -- (clu);
\draw[flow] (iso.south) -- (extract.north -| iso); \draw[flow] (clu.south) -- ++(0,-4mm) -| ($(bis.east)+(4mm,0)$) -- (bis.east);
\draw[flow] (extract) -- (out);
\draw[flow] (out) -- (final);
\node[below=5mm of final.south west, anchor=north west, align=left, font=\small, text=black] {%
  \tikz\node[box, fill=flowblue, minimum width=7mm, minimum height=3.5mm, inner sep=0pt]{};\ New preprocessing\\[1pt]
  \tikz\node[box, fill=flowgreen, minimum width=7mm, minimum height=3.5mm, inner sep=0pt]{};\ Tridiagonal MRRR code, mostly unchanged\\[1pt]
  \tikz\node[box, fill=floworange, minimum width=7mm, minimum height=3.5mm, inner sep=0pt]{};\ Modification required by Willems and Lang's theory\\[1pt]
  \tikz\node[box, fill=flowyellow, minimum width=7mm, minimum height=3.5mm, inner sep=0pt]{};\ New post-processing};
\end{tikzpicture}
\caption{High-level overview of \texttt{DBDSVDMR3}.}
\label{fig:flow}\end{figure}

\subsection{$T_{\mathrm{GK}}$ as root representation}
\label{sec:root}
On an arbitrary tridiagonal $T$, the original MRRR algorithm (i.e., Algorithm~\ref{alg:mr3}) begins by shifting to make $T$ positive definite. As explored by Willems and Lang in \cite[Section 4]{WillemsLang2012ETNA}, this shift spells disaster when applied to $T_{\text{GK}}$, effectively fusing positive and negative copies of the same singular value into one cluster. Note that this is not a problem for computing eigenvectors of $T_{\text{GK}}$, which will still come out numerically orthogonal, but \textit{is} a problem if we hope to read off left/right singular vectors according to \eqref{eqn: GK_eigenvectors}. 

\indent The remedy of Willems and Lang, adopted verbatim in \texttt{DBDSVDMR3}, is to skip this shift entirely, recalling that the entries of $B$ determine its singular values to high relative accuracy \cite{DemmelKahan1990} and therefore that $T_{\mathrm{GK}}$ is an RRR for its full spectrum. 
In terms of the code, this is accomplished by simply deleting the root-selection stage of \texttt{DSTEGR}; our implementation works directly with $T_{\mathrm{GK}}$, whose root eigenvalues $\pm\sigma_i$ are supplied by applying dqds to $B$. \\  
\indent Using $T_{\text{GK}}$ at the root level has the additional benefit of guaranteeing GK structure at the start of the representation tree. Hence, assuming that subsequent representations are computed with mixed relative stability and that element growth is limited (these are assumptions \textsc{shiftrel} and \textsc{elg} in Table~\ref{tab:conform}) the nearly-constant-diagonal condition should hold throughout. While Willems and Lang suggest some tools for monitoring element growth to avoid loss of GK-structure in \cite{WillemsLang2012ETNA}, these are left to future updates of \texttt{DBDSVDMR3}. While they were not needed to pass the LAPACK test suite, we expect these modifications to be the remedy for the accuracy losses on heavily clustered spectra reported in Section~\ref{sec:stc}, which is why they are the first item on our roadmap. \\
\indent Since MRRR struggles with tight eigenvalue clusters, it is often recommended that a small (relative) random perturbation be applied at the root level (and possibly at subsequent levels as well), which we expect to split tight clusters with high probability (see \cite[Section 5.4]{LAWN163}). We omit such a perturbation in \texttt{DBDSVDMR3} for the following reasons. First, if the perturbation is applied carelessly, specifically if we follow the recommendations of \cite{LAWN163} directly at subsequent levels, we may jeopardize the nearly constant diagonal condition, which is important here but not in the tridiagonal setting. More importantly, we found that applying small entrywise, relative perturbations  -- 
either to $T_{\mathrm{GK}}$ at the root, as \texttt{DSTEGR} does, or to the child representations of tight clusters -- did not significantly improve the results, and in fact in certain situations made them worse. 
\subsection{Full-length eigenvectors}\label{sec:vec}
\label{sec:trunc}
In MRRR for the symmetric tridiagonals, each eigenvalue yields an eigenvector whose numerically significant entries occupy a single, contiguous index window.  
\texttt{DSTEGR} takes advantage of this by computing and storing only this window -- a potentially significant optimization, lowering the work per vector from the matrix dimension to the length of the window. Nevertheless, this was found to be problematic when applied to $T_{\text{GK}}$. To be specific, the original version of \texttt{DSTEGR} stops computing an eigenvector if two consecutive entries are small enough. Since the eigenvectors of $T_{\text{GK}}$ interleave left/right singular vectors, this can prompt the algorithm to terminate too soon (e.g., if the entries of $u$ are exponentially decreasing while the entries of $v$ are exponentially increasing). We note that it is still true that the nonzero entries of the eigenvectors of a Golub-Kahan matrix are constrained to a contiguous index window -- and hence a less aggressive truncation strategy may provide some savings -- but for simplicity, and to avoid any potential failures, we disable this optimization in \texttt{DBDSVDMR3}.

\section{Numerical Results}
\label{sec:acc}
In this section, we present performance data (in terms of both accuracy and efficiency) for \texttt{DBDSVDMR3}. We compare it with four alternative solvers:\ the three existing LAPACK options for the bidiagonal SVD (see Table~\ref{tab: algorithm_comp}) and \texttt{MR3GK}, Willems and Lang's untouched tridiagonal MRRR code applied to $T_{\mathrm{GK}}$ as a black box. Our  preprocessing step is not applied to \texttt{MR3GK}; in that algorithm, $T_{\mathrm{GK}}$ is formed directly from $B$ and its main subroutine \texttt{XMR}, which does its own scaling and splitting, is asked for the $n$ largest eigenpairs, from which $u_i$ and $v_i$ are extracted by  de-interleaving and multiplying by $\sqrt{2}$. Our test suite consists of two matrix collections:\ (1) a LAPACK-style test set of 16 classes and (2) STCollection, a set of practical matrices (obtained from various applications) and stress tests maintained by Marques \cite{STCollection,MarquesDemmelVoemel2008}. Three main results emerge from these experiments:
\begin{enumerate}
\item \texttt{MR3GK} fails badly on even simple cases, which \texttt{DBDSVDMR3} avoids thanks to the implementation choices discussed in the previous section. 
\item Across the board, \texttt{DBDSVDMR3} is less accurate than the LAPACK solvers, but still provides good results on most problems. On the LAPACK test suite its median error ratios (0.03 residual, 0.18 orthogonality) match those of the LAPACK solvers, while on STCollection they are 4.2 and 6.6, about ten times the 0.3--0.7 of the LAPACK solvers. It fails (a ratio above 100) on 15 of the 125 STCollection matrices, all featuring tightly clustered eigenvalues, on which the LAPACK solvers' ratios never exceed 4.4.
\item \texttt{DBDSVDMR3} is one to two orders of magnitude faster than \texttt{DBDSVDX} on every nontrivial spectrum, and faster than \texttt{DBDSDC} wherever deflation\footnote{Deflation refers to the merge step of divide-and-conquer, where singular values of the two halves that coincide, or whose coupling coefficients are negligible, are removed from the secular equation without work. Hence, clustered or repeated spectra are cheap for \texttt{DBDSDC} while all-distinct arithmetic spectra cost the full merge.} does not help divide-and-conquer, including 25 of the 35 largest practical matrices (see Figures.~\ref{fig:pooled}--\ref{fig:tail}). At $n = 2000$, the speedups over \texttt{DBDSVDX} are $96\times$ on arithmetic spectra -- i.e., where singular values are equally spaced between $\varepsilon$ and $1$, so that nothing deflates -- and $116\times$ on the 32-decade log-uniform spectrum, while the speedups over \texttt{DBDSDC} on arithmetic spectra are $12$--$13\times$. On the aforementioned 25 large matrices, the median speedup over \texttt{DBDSDC} is $2.9\times$. These are reference-BLAS figures; when using optimized BLAS, \texttt{DBDSDC} becomes the faster full-SVD solver on all but deflation-free spectra, while \texttt{DBDSVDMR3} remains faster when only a subset of the triplets is wanted (Section~\ref{sec:blas}).
\end{enumerate}

\subsection{Experimental setup} We start by summarizing the setup for our tests.
\label{sec:setup}
\subsubsection{Notation and metrics} Throughout this section, $\mathrm{ulp}=2^{-52}$ is the double-precision machine epsilon, $\|\cdot\|_1$ is the matrix 1-norm, $I$ is the identity matrix, and \texttt{info} is the integer status returned by an LAPACK routine, which is zero unless a problem is detected. For each computed decomposition, we report three error ratios:
\begin{enumerate}
    \item Residual Error: $\|\Sigma - U^TBV\|_1/(\text{ulp} \cdot n\|B\|_1)$.
    \item Orthogonality Error: $\|I - U^TU\|_1 / (\text{ulp} \cdot n)$ and $\|I -V^TV\|_1 / (\text{ulp} \cdot n)$.
\end{enumerate}
These can be computed by LAPACK's test routines \texttt{DBDT04} and \texttt{DORT01}. When all three are on the order of one, we have accuracy at the level of floating-point roundoff. In the tables and summaries below, we count a computed decomposition as a \emph{failure} when any of the three ratios exceeds 100, equivalently when the error is more than two orders of magnitude above roundoff level. Once again, we emphasize that we do not request more than absolute accuracy from the singular values computed by \texttt{DBDSVDMR3}. 

\subsubsection{Software versions} \texttt{DBDSVDMR3} is the version merged into LAPACK (pull request 1361, merge commit \texttt{330da16} of 9 October 2026). \texttt{DBDSDC}, \texttt{DBDSQR}, and \texttt{DBDSVDX} are taken from the same source tree, the LAPACK development version that followed release 3.12.1 (git commit \texttt{51b3494} of 6 August 2026). In particular, the version of \texttt{DBDSDC} used includes the correction of a deflation-ordering bug in \texttt{DLASD2}/\texttt{DLASD7} merged in March 2026 (LAPACK pull request 1140), which is absent from release 3.12.0. \texttt{DBDSQR} and \texttt{DBDSVDX} are numerically unchanged since 3.12.0. \texttt{MR3GK} is the \texttt{XMR} code accompanying \cite{WillemsLang2012ETNA}, compiled unmodified. Unless stated otherwise, all five solvers are linked against the reference BLAS of the same source tree and run single-threaded. Section~\ref{sec:blas} repeats the timing experiments with OpenBLAS.

\subsubsection{Test suite} The first set of matrices in our test suite are constructed using LAPACK's test program \texttt{DCHKBD}. This yields 16 classes of bidiagonals:\ seven with prescribed singular-value distributions, five which are diagonal, and four that are derived from a symmetric tridiagonal matrix $T$ with a prescribed eigenvalue distribution. In the last case, we obtain a bidiagonal by shifting $T$ to be positive definite (using a Gershgorin bound) and taking its (upper) Cholesky factor. The sizes considered range over $\{2,10,30,60,100\}$ for accuracy tests and $\{100, 200, 500, 1000, 2000\}$ for efficiency tests, resulting in a total of 80  matrices in both cases.  \texttt{MR3GK} is evaluated on 75 of the 80 matrices, as it does not terminate on the zero matrix (class 1) at any size.

\indent The remaining test problems are drawn from STCollection, which contains 127 matrices (19 bidiagonals and 108 tridiagonals). These are either obtained from practical applications or are constructed as stress tests specifically for MRRR. The latter includes the glued Wilkinson matrices of \cite{DhillonParlettVoemel2006}. Once again, tridiagonal test matrices are converted to bidiagonals by shifting and computing a Cholesky decomposition. To compare accuracy, we run all five solvers on the same 122 matrices from STCollection: every one with $n\le7000$  except the three largest ($n\ge5892$), on which \texttt{DBDSQR}, whose vector computation costs $O(n^3)$, did not finish within our time budget. Results for these three test cases are still reported for the other four solvers. The remaining two matrices, with $n\approx9800$, exceed the 4\,GB host for every solver. For efficiency tests, we focus only on the 35 matrices in STCollection with sizes between 2000 and 6245. 


\subsubsection{Timing} In efficiency experiments, we time the solver call only, after one untimed warm-up call, repeating until at least $0.1$\,s has accumulated and reporting the mean; a call whose warm-up exceeds one second is a single measurement. 
The host is one unpinned virtualized Intel Xeon core (2.8\,GHz, 4\,GB) with reference BLAS and gfortran~12 at \texttt{-O2}. Absolute times drift between runs on such a host while ratios between solvers are stable, so we focus on the latter.

\subsection{Accuracy on LAPACK test matrices}
\label{sec:acc16}
\begin{figure}[!t]\centering
\includegraphics[width=\linewidth]{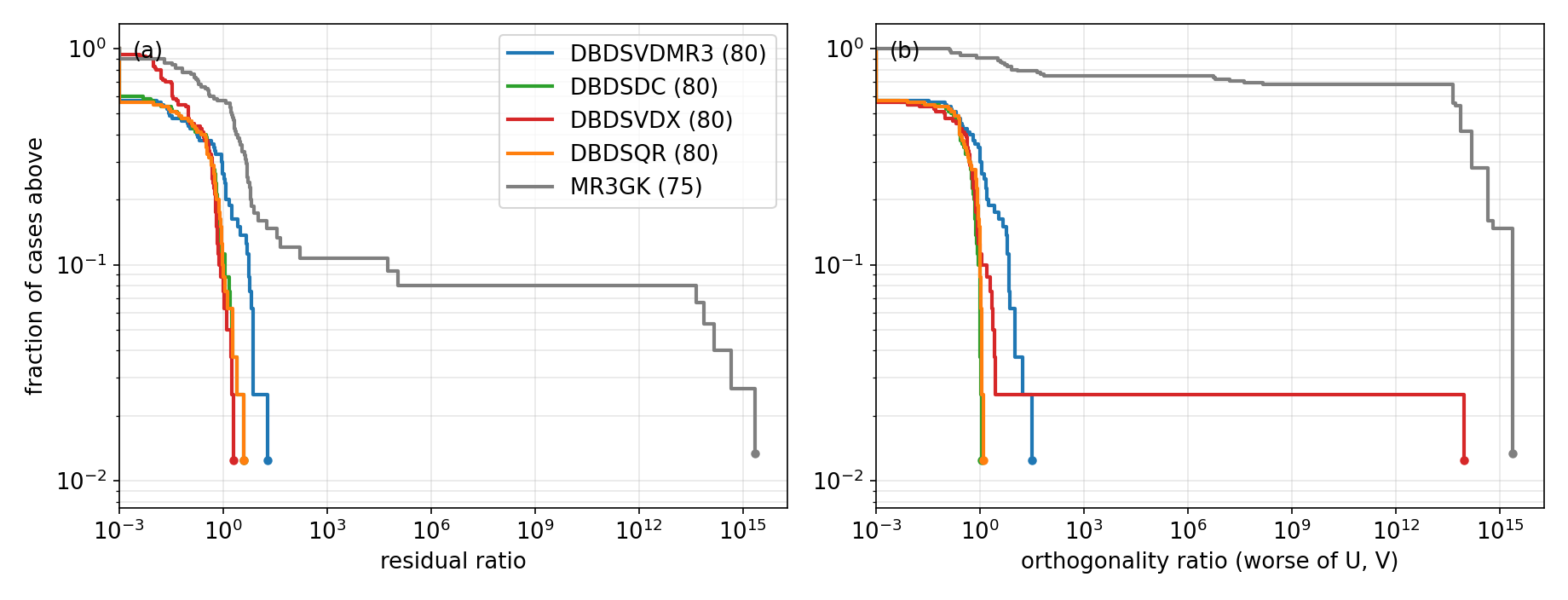}
\caption{LAPACK test suite:\ empirical distributions of (a) the residual ratio and (b) the larger of the two orthogonality ratios, one value per matrix per solver, logarithmic axes. Each curve gives the fraction of matrices whose ratio \emph{exceeds} the abscissa, so that the right end of each curve (marked by a dot, at height $1/\mbox{number of matrices}$) is that solver's worst case. \texttt{DBDSVDMR3}, \texttt{DBDSDC}, and \texttt{DBDSQR} coincide over most of the range; \texttt{MR3GK} separates from them by fifteen orders of magnitude.}
\label{fig:acccdf}\end{figure}

For the 80 test matrices constructed from \texttt{DCHKBD} and each solver, Figure~\ref{fig:acccdf} gives the full distributions of the error ratios while Table~\ref{tab:accvar} records the worst ratio per matrix variety. We note here that \texttt{DBDSDC} and \texttt{DBDSQR} stay below 4.1 (for all three error ratios) on every matrix. \texttt{DBDSVDMR3} matches them on fifteen varieties (median residual 0.03, median orthogonality 0.18); its larger values, up to 20 in residual and 31 in orthogonality, all occur on the Cholesky-derived geometric variety, where two-thirds of the singular values lie in one cluster with no relative gap above \textit{gaptol.}  Here gaptol $=\min(10^{-3},\,1/(2n_k))$ for a block of size $n_k$, the LAPACK~3.0 rule of \texttt{DLARRV} inherited from \texttt{DSTEGR}. Resolving such a cluster requires a long sequence of shifted child representations, and since \texttt{DBDSVDMR3} accepts each shift without the element-growth or nearly-constant-diagonal test that the theory of Section~\ref{sec:root} calls for, the loss of accuracy is to be expected.  \texttt{MR3GK}, the same tridiagonal kernels run as a black box, reaches orthogonality errors of $2\times10^{15}$ on every variety that contains small singular values and does not terminate on the zero matrix.
\texttt{DBDSVDX} has one large value -- an orthogonality error of $9\times10^{13}$ on the variety with singular values spread log-uniformly over 32 decades. Half of those singular values lie below $\mathrm{ulp}\cdot\|B\|$, closer to one another and to their negatives than bisection with an absolute tolerance can resolve, so inverse iteration returns non-orthogonal vectors for them.


\begin{table}[!t]\caption{LAPACK test suite:\ worst ratio per matrix variety, maximized over the five sizes $\left\{ 2, 10, 30, 60, 100 \right\}$. For \texttt{DBDSVDMR3} the residual, $U$-orthogonality, and $V$-orthogonality ratios are shown separately; for the other solvers we record only the largest of the three. Variety 7 is the Cholesky-derived geometric spectrum discussed in the text.}
\label{tab:accvar}\centering\footnotesize\setlength{\tabcolsep}{3pt}
\begin{tabular}{rl|c|rrrr}
\toprule
 \multicolumn{2}{l|}{Variety} & \texttt{DBDSVDMR3} (res / $U$ / $V$) & \texttt{DBDSDC} & \texttt{DBDSVDX} & \texttt{DBDSQR} & \texttt{MR3GK} \\
\midrule
1 & zero & 0 / 0 / 0 & 0 & 0 & 0 & hang \\
2 & identity & 0 / 0 / 0 & 0 & 0.5 & 0 & 4 \\
3 & arith bidiag & 6.53 / 6.82 / 6.78 & 1 & 0.736 & 1.34 & 2.25e15 \\
4 & geom bidiag & 1.23 / 1 / 1.12 & 1.78 & 0.697 & 1 & 2.25e15 \\
5 & log-dist & 0.897 / 1 / 1 & 1 & 9.4e13 & 1 & 2.25e15 \\
6 & arith tri-chol & 7.53 / 10.1 / 10 & 1.85 & 2.57 & 1.94 & 8.32 \\
7 & geom tri-chol & 19.8 / 30.8 / 30.8 & 1.53 & 2.81 & 0.599 & 68.4 \\
8 & clust tri-chol & 0.578 / 0.292 / 0.625 & 0.15 & 0.643 & 0.15 & 5.25 \\
9 & arith sqrt-ovfl & 5.64 / 6.26 / 6.85 & 1.14 & 1.12 & 1.14 & 2.25e15 \\
10 & arith sqrt-unfl & 5.44 / 6.13 / 6.04 & 2.47 & 1.06 & 2.47 & 2.25e15 \\
11 & arith diag & 0 / 0 / 0 & 0 & 0.5 & 0 & 2.25e15 \\
12 & geom diag & 0 / 0 / 0 & 0 & 0.5 & 0 & 2.25e15 \\
13 & clust diag & 0 / 0 / 0 & 0 & 0.5 & 0 & 2.25e15 \\
14 & diag sqrt-ovfl & 0 / 0 / 0 & 0.0167 & 0.25 & 0 & 2.25e15 \\
15 & diag sqrt-unfl & 0 / 0 / 0 & 0 & 0.5 & 0 & 2.25e15 \\
16 & SPD geom-chol & 1.19 / 1 / 1.12 & 4.03 & 0.756 & 4.03 & 1.45e08 \\
\bottomrule
\end{tabular}
\end{table}

\subsection{Accuracy on STCollection}
\label{sec:stc}
\begin{figure}[!t]\centering
\includegraphics[width=\linewidth]{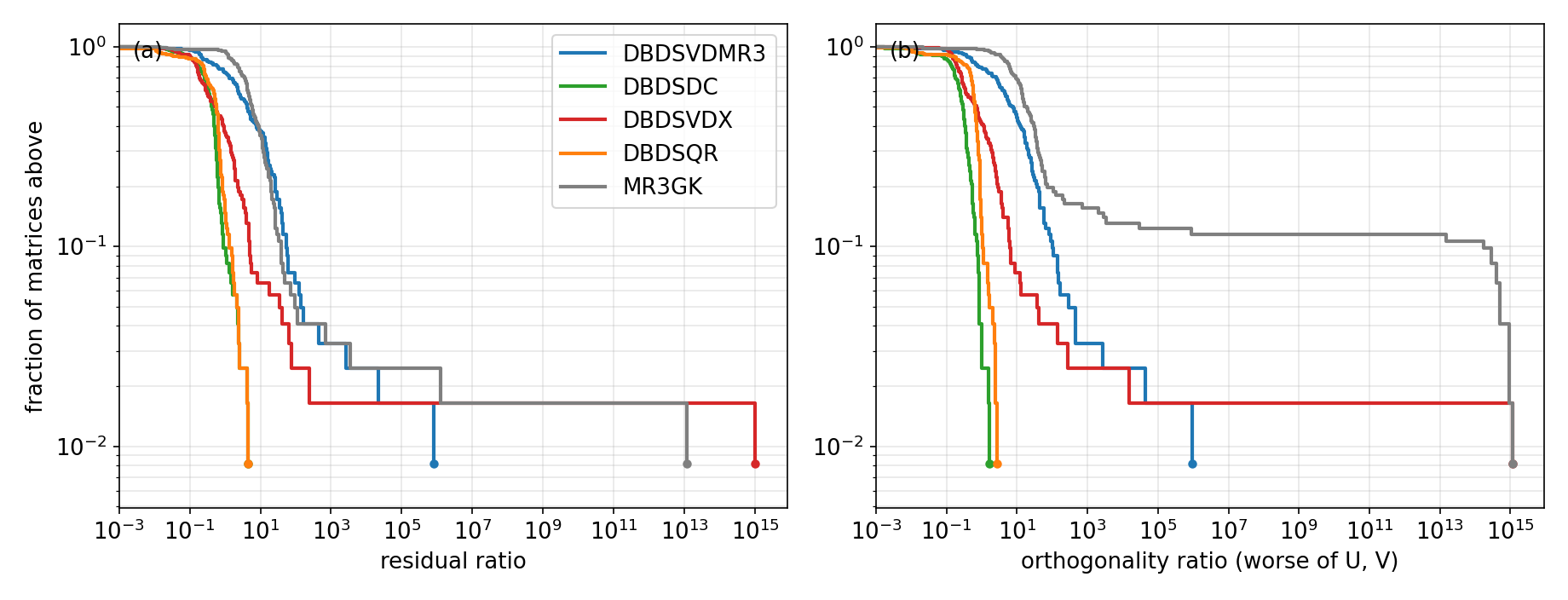}
\caption{Practical suite, the same 122 matrices for all five solvers ($n$ up to 5472):\ empirical distributions of (a) the residual ratio and (b) the larger of the two orthogonality ratios, one value per matrix per solver, logarithmic axes, drawn as in Figure~\ref{fig:acccdf}. That is, each curve gives the fraction of matrices whose ratio exceeds the abscissa and ends at the solver's worst case. \texttt{DBDSVDMR3}'s curves lie to the right of the LAPACK solvers' throughout; its far right end is the glued Wilkinson family, while the far right end of \texttt{MR3GK}'s is the 10\% of matrices on which black-box MR$^3$ loses orthogonality completely. \texttt{DBDSVDX}'s two isolated large values are \texttt{B\_bug414} (underflow, $\texttt{info}=2n+1$) and \texttt{B\_16}; \texttt{DBDSDC} and \texttt{DBDSQR} have none.}
\label{fig:stress}\end{figure}

Figure~\ref{fig:stress} compares the five solvers on the 122 matrices from STCollection mentioned above. \texttt{DBDSQR} achieves no error ratio above 4.4 on any matrix. In most cases, the accuracy ordering is as follows: (1) \texttt{DBDSDC} (median residual 0.41, median orthogonality 0.29), (2) \texttt{DBDSVDX} and \texttt{DBDSQR} (0.5--0.7 on both), (3) \texttt{DBDSVDMR3} (4.2 and 6.6), and (4) \texttt{MR3GK} (5.7 and 16). For practical matrices, the necessity of the root modification shows in the tail of the distributions -- where \texttt{MR3GK} exceeds $10^{13}$ on 10\% of the matrices while \texttt{DBDSVDMR3} never reaches this threshold --- rather than in the average case. \\
\indent \texttt{DBDSVDMR3}'s largest ratios all arise on heavily clustered spectra, with by far the largest occurring on glued Wilkinson matrices, in which the 21 eigenvalues of $W_{21}^{+}$, the $21\times21$ Wilkinson matrix, are repeated a hundred times with a coupling between $10^{-14}$ and $10^{14}$. In this case, residual/orthogonality ratios reach $8\times10^{5}$ and $9\times10^{5}$ (\texttt{T\_W21\_g\_1e+02}, $n{=}2100$), while on the same matrix \texttt{DBDSDC} and \texttt{DBDSQR} stay below 3, \texttt{DBDSVDX} below 7, and even \texttt{MR3GK} below 5. This failure is \texttt{DBDSVDMR3}'s alone. On the three largest matrices, not covered by Figure~\ref{fig:stress}, the same picture holds (\texttt{T\_bcsstkm13\_3}, $n{=}6009$: \texttt{DBDSVDMR3} $3.2\times10^{3}$ and $4.4\times10^{3}$, \texttt{DBDSDC} 0.2 and 0.2, \texttt{DBDSVDX} 3.2 and 4.2). These matrices are the classic MRRR stress case \cite{DhillonParlettVoemel2006}, precisely because their clusters force the deepest sequences of shifts. With those shifts unguarded, \texttt{DBDSVDMR3} loses accuracy in proportion to the depth of the representation tree. While singular value clustering is where its errors live, it is not by itself the cause; many equally clustered matrices in the collection are computed at reference quality, and identifying which cluster descents lose accuracy requires the instrumented shift-acceptance test that we have not yet implemented. This is the solver's main limitation in its present form, compounded by a second:\ it returns $\texttt{info}=0$ on every one of these matrices, so the caller is not warned.

Of the other solvers, \texttt{DBDSDC} and \texttt{DBDSQR} fail on no matrix of either suite. Both return $\texttt{info}=0$ throughout, with no error ratio above 4.4. We note that the \texttt{DBDSDC} of LAPACK release 3.12.0 does fail, returning $\texttt{info}=1$ (the divide-and-conquer update did not converge) and meaningless failures on \texttt{B\_Kimura\_429} and \texttt{T\_nasa1824\_3}; the correction to \texttt{DLASD2}/\texttt{DLASD7} in LAPACK pull request 1140 (March 2026), included in the version used here, resolves both. \texttt{DBDSVDX} is the one LAPACK solver with failure modes in our tests.  It performs no scaling, and consequently on \texttt{B\_bug414}, whose entries are as small as $10^{-171}$, its squared quantities underflow and it stops with an internal error (\texttt{info}$=2n+1$). Further, on \texttt{B\_16}, whose singular values span sixty decades, and on the log-uniform class from the LAPACK set, the singular values below $\mathrm{ulp}\cdot\|B\|$ cannot be separated by bisection with an absolute tolerance, and hence the vectors computed by inverse-iteration come out non-orthogonal (orthogonality ratio $1.5\times10^{4}$). Finally, \texttt{MR3GK} behaves as poorly here as on the LAPACK test suite.


\subsection{Timing}
\label{sec:perf}
\begin{figure}[!t]\centering
\includegraphics[width=.85\linewidth]{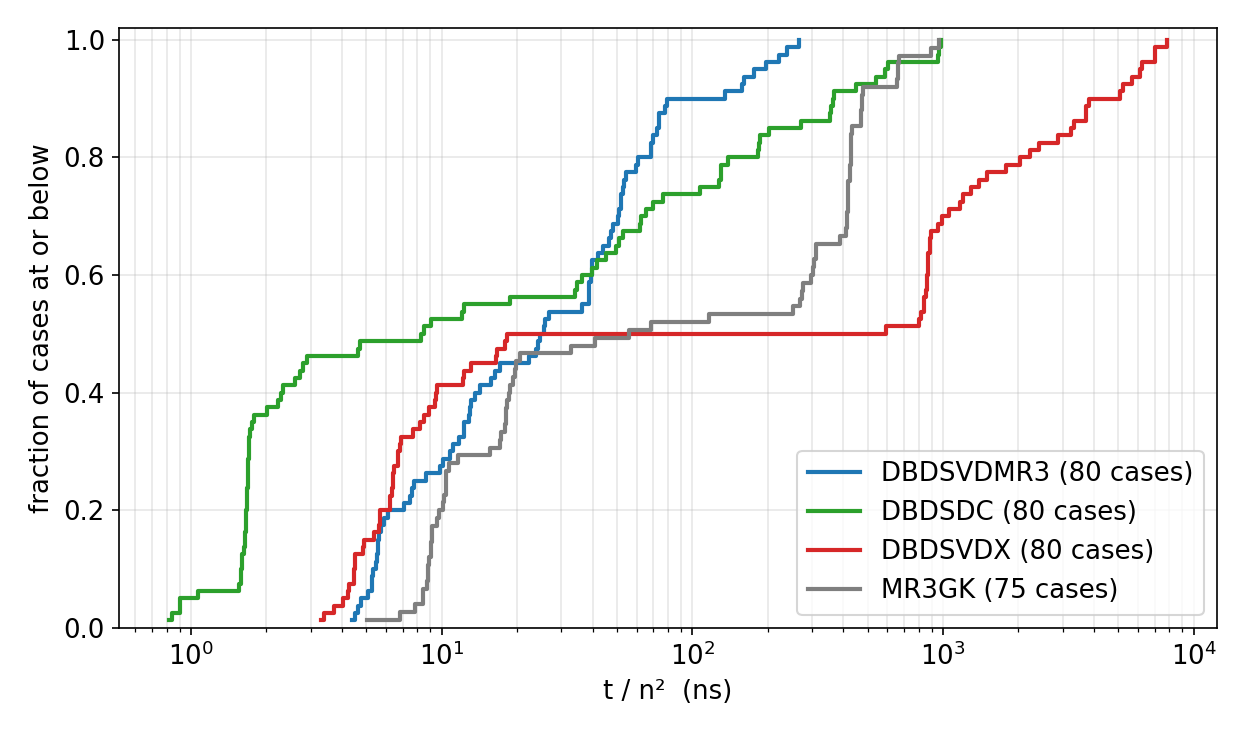}
\caption{LAPACK test suite: empirical distribution of the normalized solve time $t/n^2$ over the 80 (variety, size) cases, $n\in\{100,\dots,2000\}$ (\texttt{MR3GK}: 75). \texttt{DBDSDC} has the lowest median, \texttt{DBDSVDMR3} the lowest maximum, and \texttt{DBDSVDX} lies an order of magnitude to the right.}
\label{fig:pooled}\end{figure}

Figure~\ref{fig:pooled} shows the distribution of normalized time $t/n^2$ over the 80 cases in the LAPACK test suite. Against \texttt{DBDSVDX}, \texttt{DBDSVDMR3} is faster by one to two orders of magnitude on every class with a nontrivial spectrum ($96\times$ on the arithmetic variety at $n{=}2000$, $116\times$ on the 32-decade log-uniform one). \texttt{DBDSVDX} is faster only on the zero, identity, and diagonal varieties, which it dispatches by early exit, and on the clustered Cholesky-derived variety. Against \texttt{DBDSDC}, the outcome depends on the spectrum. \texttt{DBDSVDMR3} is $12$--$13\times$ faster at $n{=}2000$ on the arithmetic varieties, where nothing deflates, and $2.8\times$ and $1.4\times$ faster on the Cholesky-derived arithmetic and SPD varieties, but it is also $1.0$--$4.1\times$ slower on the geometric and clustered varieties and $8$--$15\times$ slower on the zero, identity, and diagonal varieties. 
Over all LAPACK-constructed tests, \texttt{DBDSVDMR3} has the smallest maximum of $t/n^2$ ($266$ against $981$\,ns) and the smaller total time at $n{=}2000$ ($4.2$\,s against $14.6$\,s), while \texttt{DBDSDC} has the smaller median.

\begin{figure}[!t]\centering
\includegraphics[width=.85\linewidth]{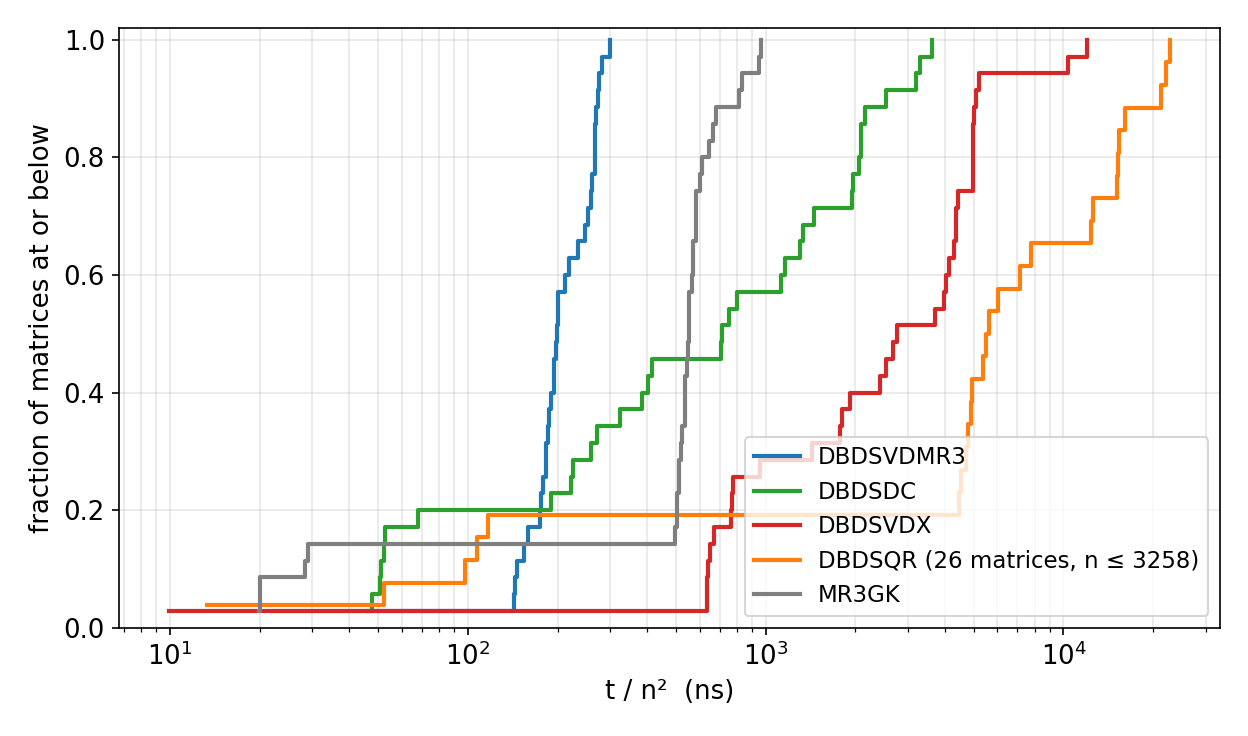}
\caption{Practical suite, the 35 matrices with $n\ge2000$: empirical distribution of the normalized solve time $t/n^2$ per matrix, single measurements. \texttt{DBDSVDMR3}, \texttt{DBDSDC}, \texttt{MR3GK}, and \texttt{DBDSVDX} were timed on all 35 ($n$ up to 6245, at which point \texttt{DBDSVDX} takes 404\,s). \texttt{DBDSQR} is shown on the 26 matrices with $n\le3258$; its vector computation grows like $n^3$ and a single run takes 235\,s at $n=3258$, 14 minutes at $n=4704$, and 26 minutes at $n=5472$, beyond our budget for the largest matrices. Medians of $t/n^2$: 199\,ns (\texttt{DBDSVDMR3}), 551 (\texttt{MR3GK}), 711 (\texttt{DBDSDC}), 2767 (\texttt{DBDSVDX}), 5567 (\texttt{DBDSQR}).}
\label{fig:stctime}\end{figure}

\begin{figure}[!t]\centering
\includegraphics[width=\linewidth]{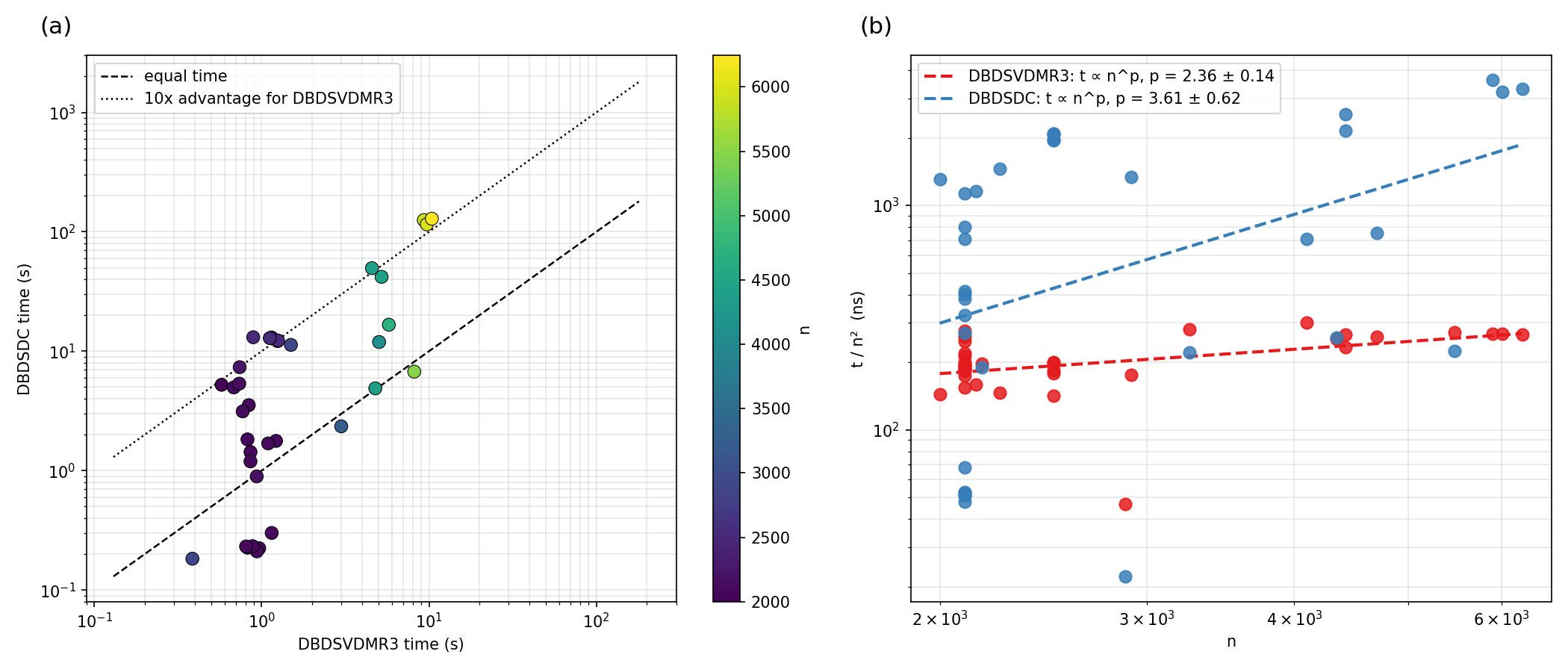}
\caption{Practical suite, the 35 matrices with $n\ge2000$. (a) Solve time of \texttt{DBDSDC} (vertical) against that of \texttt{DBDSVDMR3} (horizontal), one point per matrix colored by $n$, logarithmic axes; the dashed line marks equal time and the dotted line a tenfold advantage for \texttt{DBDSVDMR3}. Twenty-five points lie above the diagonal, the largest matrices among them by an order of magnitude. (b) The same runs as $t/n^2$ against $n$, with a power law $t\propto n^{p}$ fitted to each solver ($p$ and its standard error in the legend). \texttt{DBDSVDMR3}'s normalized time is nearly constant ($p=2.36\pm0.14$), whereas \texttt{DBDSDC}'s grows with $n$; its exponent is poorly determined ($3.61\pm0.62$) because its cost is governed by deflation rather than by $n$.}

\label{fig:tail}\end{figure}

Two comparisons are made on the practical matrices in STCollection with $n\ge2000$. Figure~\ref{fig:stctime} shows each solver over the matrices it was timed on; all solvers except \texttt{DBDSQR} were timed on all 35, and \texttt{DBDSQR}, whose cost grows like $n^3$ (26 minutes for a single run at $n{=}5472$), on the 26 with $n\le3258$. \texttt{MR3GK} is slower than \texttt{DBDSVDMR3} on 30 of the 35 matrices (median factor $2.5$). The ordering by median normalized time is \texttt{DBDSVDMR3} (199\,ns), \texttt{MR3GK} (551), \texttt{DBDSDC} (711), \texttt{DBDSVDX} (2767), \texttt{DBDSQR} (5567). Figure~\ref{fig:tail} compares \texttt{DBDSVDMR3} and \texttt{DBDSDC} alone on all 35 matrices ($n$ up to 6245), where the other solvers were not affordable. In this setting, \texttt{DBDSVDMR3} is faster on 25 of the 35, by a median factor of $2.9$ and up to $14.7$, and the totals are 89\,s against 618\,s. The largest margins are on structural-engineering matrices -- e.g., \texttt{bcsstkm11\_4} ($n{=}5892$), 9.3\,s against 125.3\,s; \texttt{bcsstkm13\_3} ($n{=}6009$), 9.7\,s against 115.5\,s; and \texttt{Alemdar\_1} ($n{=}6245$), 10.4\,s against 128.7\,s -- and the ten matrices on which \texttt{DBDSDC} is faster (by up to $4\times$ on the glued Wilkinson matrices, $2\times$ on \texttt{zenios}) are those whose spectra deflate well. 

\subsubsection{Optimized BLAS and subsets of the spectrum}\label{sec:blas}
\begin{figure}[!t]\centering
\includegraphics[width=\linewidth]{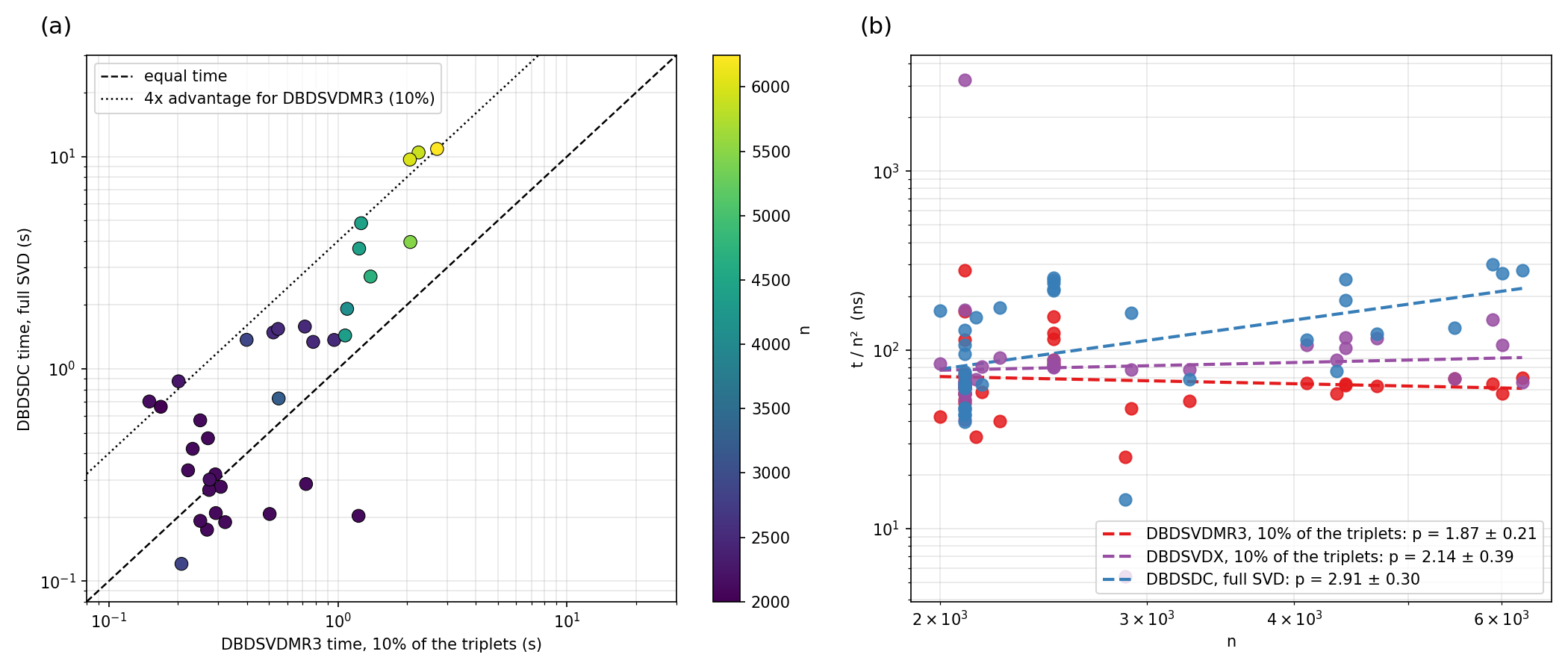}
\caption{As Figure~\ref{fig:tail}, with single-threaded OpenBLAS; \texttt{DBDSVDMR3} (follow-up version) and \texttt{DBDSVDX} compute only the largest $10\%$ of the triplets; \texttt{DBDSDC}, which has no subset option, computes all of them. (a) 25 of 35 points lie above the equal-time line, i.e., \texttt{DBDSVDMR3} is faster. (b) Fitted exponents $1.87\pm0.21$ (\texttt{DBDSVDMR3}), $2.14\pm0.39$ (\texttt{DBDSVDX}), $2.91\pm0.30$ (\texttt{DBDSDC}).}
\label{fig:subset}\end{figure}

All timings above use the reference BLAS. With single-threaded OpenBLAS~0.3.26 (LAPACK from the same source tree, on a second host where a reference-BLAS control reproduced the ratios of Figures~\ref{fig:pooled}--\ref{fig:tail}), \texttt{DBDSVDMR3}, whose kernels are BLAS-1, runs at the same speed, while \texttt{DBDSDC}, whose work is in matrix--matrix products, becomes five times faster. For the full SVD, \texttt{DBDSDC} is then faster on 31 of the 35 STCollection matrices with $n\ge2000$, by $5$--$10\times$ where deflation is heavy and by under $20\%$ on the three largest problems. On spectra without deflation, \texttt{DBDSVDMR3} stays ahead by a margin that grows with $n$ (see Table~\ref{tab:deflfree}). \\
\indent This matches what is known for the tridiagonal eigenproblem:\ with optimized BLAS, divide-and-conquer is expected to win if the full spectrum is desired and deflation is substantial, while MRRR wins for spectra without deflation or, because of its $O(nk)$ complexity, when only a subset of the spectrum is requested \cite{DemmelMarquesParlettVoemel2008,DhillonParlettVoemel2006,PetschowBientinesi2011}. In the latter case, the comparison is uneven by design. That is, \texttt{DBDSDC} and \texttt{DBDSQR} always compute the full SVD, while \texttt{DBDSVDX} accepts a \texttt{RANGE} argument. The merged \texttt{DBDSVDMR3} also computes all vectors regardless of \texttt{RANGE}, though in a follow-up pull request we intend to make an adaptive version (that only computes the requested triplets) available. Figure~\ref{fig:subset} shows the result when only the largest $10\%$ of the triplets are requested. In this case, computing the subset costs $4.7\times$ less than the full SVD and \texttt{DBDSVDMR3} is faster than \texttt{DBDSDC} on 25 of the 35 matrices (by $4$--$4.7\times$ on the large structural-engineering matrices, with \texttt{DBDSDC} winning only on the glued-Wilkinson family, which it deflates to a full SVD in $0.2$\,s) and is about as fast as \texttt{DBDSVDX} on the same request (median time ratio 1.06, totals 26\,s against 48\,s). The subset path also makes \texttt{DBDSVDMR3} embarrassingly parallel in the manner of MRRR for the tridiagonal eigenproblem \cite{BientinesiDhillonVanDeGeijn2005,Voemel2010}, where each of $p$ processors calls it with \texttt{RANGE='I'} on a disjoint block of about $n/p$ indices. This requires no communication, with the only redundant work 
corresponding to representations of clusters that straddle a block boundary. 

\begin{table}[!t]\centering\footnotesize
\caption{Full SVD with single-threaded OpenBLAS on two problems without deflation:\ an equally spaced spectrum (bidiagonal from \texttt{DLATMS}) and the Toeplitz bidiagonal $d_i{=}e_i{=}1$. In both cases, ratio = \texttt{DBDSDC} time / \texttt{DBDSVDMR3} time, so values above one favor \texttt{DBDSVDMR3}. On random, geometric, and clustered spectra the ratio is below one.}
\label{tab:deflfree}
\begin{tabular}{llrrr}
\toprule
spectrum & $n$ & \texttt{DBDSVDMR3} & \texttt{DBDSDC} & ratio \\
\midrule
equally spaced, $\kappa=1/\mathrm{ulp}$ & 2000 & 0.29\,s & 0.54\,s & 1.8 \\
 & 4000 & 1.25\,s & 2.98\,s & 2.4 \\
 & 8000 & 6.59\,s & 18.85\,s & 2.9 \\
Toeplitz $d_i{=}e_i{=}1$ & 2000 & 0.44\,s & 0.76\,s & 1.8 \\
 & 4000 & 1.59\,s & 4.86\,s & 3.1 \\
 & 8000 & 7.67\,s & 32.51\,s & 4.2 \\
\bottomrule
\end{tabular}
\end{table}

In total, these results answer our main question with unusual precision:\ MRRR can, in fact, be used (with minimal modification) to compute the bidiagonal SVD via the Golub-Kahan matrix. Indeed, \texttt{DBDSVDMR3} passes the LAPACK test suite with no error ratio above 31 and keeps all three ratios below 100 on 88\% of STCollection, while outrunning \texttt{DBDSVDX} and \texttt{DBDSQR} by one to two orders of magnitude. With the reference BLAS it also beats divide-and-conquer, the current best option for large matrices, in both worst-case time and total time over the suite. With optimized BLAS it keeps that advantage (1) on spectra that do not deflate and (2) whenever only part of the spectrum is requested (once its subset path is merged).  


\section{Recommendations for AI-assisted code development}
\label{sec:ai}
The code modifications described in Section~\ref{sec:adapt} were produced with an AI coding assistant (Claude Code \cite{ClaudeCode}) working under the direction of a numerical linear algebra expert, who did not write the Fortran but fixed the mathematical constraints through prompts. Two instructions did most of the work:\ (1) that the Golub--Kahan matrix is a distinct problem from the general symmetric tridiagonal one, so that the positive definite root shift is to be removed rather than tuned and (2) that the code should stay as faithful to the original MR$^3$ implementation as possible, changing nothing the theory does not require. Each proposed change was accepted only after the accuracy and timing experiments of Section~\ref{sec:acc} had been rerun on the same input matrices, and rejected otherwise.  \\
\indent Concretely, the division of labor was as follows. The human supplied the starting point (\texttt{DSTEGR}'s kernels rather than a fresh implementation), the two instructions above, the test matrices, and the acceptance criterion. The assistant then wrote all of the Fortran, the test harness, and the extensions of LAPACK's test suite. Of the modifications in Section~\ref{sec:adapt}, the use of $T_{\mathrm{GK}}$ as the root representation (Section~\ref{sec:root}) was prescribed by the human as a mathematical requirement, while the preprocessing and splitting of Section~\ref{sec:prep}, the handling of odd-dimensional blocks and their null vectors, and the removal of the support-window truncation of Section~\ref{sec:vec} were proposed by the assistant when asked to diagnose and remove failures on specific test matrices. These modifications were accepted after the human had checked them against the theory of Section~\ref{sec:bg}. The random perturbation discussed in Section~\ref{sec:root} is an example of a change proposed by the assistant and later rejected on the evidence. The same workflow can be carried further --- for example to improve the quality of individual kernels such as dot-product accumulations, or to instrument the child-acceptance test discussed in Section~\ref{sec:root} --- and we regard this work as a case study in AI-assisted modification of established numerical software.\\
\indent  Finally, we emphasize that  part of the success of our effort can be attributed to the fact that one of the authors on this paper wrote the original \texttt{DSTEGR} code. Accordingly, when the AI agent suggested a technical change like the one discussed in Section~\ref{sec:trunc}, we felt confident in our understanding of both what the original code was doing and why it could be safely disabled. This suggests that AI-accelerated coding yields the best results when working with existing implementations that are well understood. Indeed, it was for this reason that we chose to build \texttt{DBDSVDMR3} on top of \texttt{DSTEGR} rather than attempting to fix the Willems and Lang code or working with \texttt{DSTEMR} instead.

\section{Availability and reproduction}
The solver, its single-precision version, its LAPACKE bindings, and the dense drivers \texttt{[SDCZ]GESVDMR3} built on it are part of the reference LAPACK repository (merged via pull request 1361). The subset path discussed in Section~\ref{sec:blas} is the subject of a follow-up pull request. The test harness, matrix lists, raw logs, and the scripts that regenerate every figure and table of this paper are in the \texttt{xmr-evolve} repository, pinned by commit.


\section*{Acknowledgments}
The solver was developed and analyzed with substantial assistance
from large-language-model agents (Anthropic Claude) operating in
evaluator-gated loops -- i.e., all AI-produced changes were verified by the
diffs, ablations, and measurements reported here, and this manuscript
was drafted with the same assistance. Per SIAM policy, AI tools are not
authors; their use is disclosed here.

\bibliographystyle{siamplain}
\bibliography{refs}
\end{document}